\documentclass[11pt]{article}
\usepackage{amssymb,latexsym}
\usepackage{url}

 \catcode`@=11\def\@oddhead{\hbox{}\hfil\rm\thepage}\def\@oddfoot{}
 \def\@evenhead{\hbox{}\hfil\rm\thepage}\def\@evenfoot{}
 \@twosidetrue \catcode`@=12
\newtheorem{prp}{Proposition}
\newtheorem{lem}[prp]{Lemma}\newtheorem{thm}[prp]{Theorem}
\newenvironment{prf}{\begin{trivlist}\item[\emph{Proof.}]}{\end{trivlist}
  \medskip\par}
\newenvironment{rem}{\begin{trivlist}\item[\emph{Remarks.}]}{\end{trivlist}
  \medskip\par}
\def\prpb{\begin{prp}}\def\prpe{\end{prp}}
\def\lemb{\begin{lem}}\def\leme{\end{lem}}
\def\thmb{\begin{thm}}\def\thme{\end{thm}}
\def\corb{\begin{cor}}\def\core{\end{cor}}
\def\prfb{\begin{prf}}\def\prfe{\end{prf}}
\def\remb{\begin{rem}}\def\reme{\end{rem}}
\def\prpa#1{\label{p:#1}}

\def\thma#1{\label{t:#1}}\def\thmu#1{Theorem~\ref{t:#1}}

\def\seca#1{\label{s:#1}}\def\secu#1{Section~\ref{s:#1}}

\def\itmb{\begin{enumerate}}\def\itme{\end{enumerate}}
\def\itdb{\begin{itemize}}\def\itde{\end{itemize}}
\def\ittb{\begin{description}}\def\itte{\end{description}}
\def\eqnb{\begin{equation}}\def\eqne{\end{equation}}
\def\arrb#1{\begin{array}{#1}}\def\arre{\end{array}}
\def\tabb#1{\par\noindent\begin{tabular}{#1}}
\def\tabe{\end{tabular}\par\noindent}
\def\eqna#1{\label{e:#1}}\def\eqnu#1{(\ref{e:#1})}

\def\QED{\relax\ifmmode\let\@tempa\relax\ifcase\@eqcnt\def\@tempa{& & &}\or
  \def\@tempa{& &}\else\def\@tempa{&}\fi\@tempa $\Box$ \else\hfill $\Box$ \fi}
\def\DDD{\relax\ifmmode\let\@tempa\relax\ifcase\@eqcnt\def\@tempa{& & &}\or
 \def\@tempa{& &}\else\def\@tempa{&}\fi\@tempa $\Diamond$
 \else\hfill $\Diamond$ \fi}

\def\Rom#1{\uppercase\expandafter{\romannumeral#1}}
\def\dsp{\displaystyle}

\def\limf#1{\displaystyle \lim_{#1\to\infty}}

\def\Ccomb#1#2{\setbox0=\hbox{$\displaystyle\mathrm{C}$}\setbox1=\hbox{%
$\scriptstyle #1$}\kern \wd1{\mathrm{C}}_{\kern -1.05\wd0\kern -0.99\wd1{#1}
 \kern 1.15\wd0{#2}}}

\def\clvec#1#2#3{\def\clvecone{#3}\left(\arrb{c} \dsp #1\\ \dsp #2
 \ifx\clvecone\empty\else\\ \dsp #3\fi\arre\right)}

\def\diff#1#2{\dsp\frac{d\,#1}{d#2}}

\def\pderiv#1#2{\dsp\frac{\partial\,#1}{\partial#2}}

\def\le{\leqq} \def\leq{\leqq}\def\ge{\geqq} \def\geq{\geqq}

\def\reals{{\mathbb R}}
\def\preals{{\mathbb R_+}}

\def\prb#1{\def\prbone{#1}
  \ifx\prbone\empty{\mathrm{P}}\else{\mathrm{P[\;}}#1{\mathrm{\;]}}\fi}
\def\prbseq#1#2{\def\prbseqone{#2}
  \ifx\prbseqone\empty{\mathrm{P}}_{#1}\ignorespaces
  \else{\mathrm{P}}_{#1}{\mathrm{[\;}}#2{\mathrm{\;]}}\fi}
\def\EE#1{{\mathrm{E[\;}}#1{\mathrm{\;]}}}
\def\EEseq#1#2{\def\EEseqone{#2}
  \ifx\EEseqone\empty{\mathrm{E}}_{#1}\else
 {\mathrm{E}}_{#1}{\dsp\mathrm{[\;}}#2{\mathrm{\;]}}\fi}

\def\VVseq#1#2{\def\VVseqone{#2}
  \ifx\VVseqone\empty{\matrm{V}}_{#1}\else
 {\mathrm{V}}_{#1}{\dsp\mathrm{[\;}}#2{\mathrm{\;]}}\fi}

\def\ssN{^{(N)}}

\title{
 Hydrodynamic limit of move-to-front rules and search cost probabilities
}
\author{
Kumiko Hattori
\\ 
\small Department of Mathematics and Information Sciences,
\\
\small  Tokyo Metropolitan University, Hachioji, Tokyo 192-0397, Japan.
\\ \small email: \url{khattori@tmu.ac.jp}
\\ \and
Tetsuya Hattori
\\ 
\small Laboratory of Mathematics, Faculty of Economics, Keio University, 
\\
\small 4--1--1 Hiyoshi, Yokohama 223-8521, Japan
\\ \small URL: \url{http://web.econ.keio.ac.jp/staff/hattori/research.htm}
\\ \small email: \url{hattori@econ.keio.ac.jp}
} 
\date{2009/07/20}

\begin{document}
\maketitle

\begin{center}
ABSTRACT
\end{center}


We study a hydrodynamic limit approach to move-to-front rules,
namely, a scaling limit as the number of items tends to infinity,
of the joint distribution of jump rate and position of items.
As an application of the limit formula, 
we present asymptotic formulas on search cost probability
distributions, applicable for general jump rate distributions.

\vspace*{1in}\par
\bigskip\par
\footnotetext{ \noindent\textit{Corresponding author:} 
Tetsuya Hattori, 
\url{hattori@econ.keio.ac.jp}
\par\noindent
Laboratory of Mathematics, Faculty of Economics, Keio University, 
4--1--1 Hiyoshi, Kohoku-ku, Yokohama 223--8521, Japan.
(Tel: \url{011-81-45-566-1300})
}

\newpage

\section{Introduction.}
\seca{1}

The {\it move-to-front} (MTF) rule is an algorithm for a self-organizing 
linear list of 
a finite number of items, say, $\{1, 2, \ldots , N\}$.  
The list is updated in the following way.  
At each discrete unit of time, an 
item is requested, according to request probability $p_i>0$, 
$i=1, \ldots , N$.  
If the item is found at the $k$th position, it is moved to the top position 
and 
items in the first to the $(k-1)$th positions are moved down by one position. 
Successive requests are independent.  This algorithm defines a Markov chain 
on the state space of the permutations of $\{1, 2, \ldots , N\}$.  
There have been extensive studies on the MTF model,   
dating back to \cite{Tsetlin1963,mv2frnt1,mv2frnt2}.

In \cite{HH071,HH072,HH073} we studied a continuous time Markov process
which we called the {\it stochastic ranking process}.
The process corresponds to a Poisson embedding of 
the MTF chain into continuous-time \cite{mv2frnt4,BlomHolst91}.  
Each item makes jumps to the 
top with jump rate per unit time $w_i$ (corresponding to 
$p_i$ in the discrete-time model) independently of the others.

Near the top of the list, popular (often-jumped or often-requested) items 
tend to gather, but 
there are always unpopular items mixed with popular ones.  
As a mathematically precise formulation of such an observation,
we proved in \cite{HH071} that,
under appropriate conditions
such as the existence of the limit jump-rate 
distribution $\lambda $ as $N \to \infty$, 
the joint distribution $\mu\ssN_t$ of the jump rate (popularity) 
and the scaled position on the list converges as $N \to \infty$,
and also gave an explicit formula for the limit distribution $\mu_t$\,.
We also obtained the expression for the boundary on the list between items 
that have jumped at least once and those that have not.   
Under an appropriate scaling, 
the boundary converges to a deterministic trajectory $y=y_C(t)$ 
as $N \to \infty$.  $y_C(t)$ is given by the Laplace transform of the 
limit jump-rate distribution $\lambda$:
\[
y_C(t)=1-\int_0^{\infty} e^{-wt} \lambda(dw).
\]
$\mu_t$ mentioned above has a general expression in terms of 
the inverse function $t_0(y)$ of $y_C(t)$ and its likes
(see \eqnu{Tets20070726integrated} or \eqnu{spacialdistributionfcn}
in \secu{2}).

After \cite{HH071,HH072} were accepted for publication, we learned that
the MTF rule has been in the literature
for nearly half a century \cite{Tsetlin1963,mv2frnt1,mv2frnt2,mv2frnt3},
and has also been called self-organizing search, 
Tsetlin library \cite{Rivest1976}, or more recently,
least-recently-used (LRU) caching \cite{Jelenkovic03,SM2006}.
In spite of a long history of studies in the rule,
the main results in \cite{HH071,HH072},
which we summarize in \secu{2}, have escaped being noticed.
Mathematical reasons why the curve $y=y_C(t)$ plays an important role
in the formula for $\mu_t$ and also why its inverse function $t=t_0(y)$
appears in $\mu_t$
(see \eqnu{Tets20070726integrated} or \eqnu{spacialdistributionfcn})
are studied in \cite{HH072}, where it is proved that
\eqnu{Tets20070726integrated} satisfies
a system of non-linear Burgers type partial differential equations (PDE),
which can be interpreted as a motion of mixed incompressible
fluid driven by evaporation.
An initial value problem for the PDE is solved 
by a standard method of characteristic curves, 
one of which is exactly the curve $y=y_C(t)$. 
The solution to the PDE is then written using
the inverse function of the characteristic curves.
In view of this result, \thmu{HDL} could be viewed as a 
mathematical result on a {\it hydrodynamic limit}.

Our formula also has a direct practical application on the web.
We noted in \cite{HH072,HH073} that
the characteristic curve $y=y_C(t)$ is 
actually observed on the internet as the time-development of web rankings,
which have become popular in the late twentieth century,
as a result of the advance in web technology.
In \cite{HH072,HH073} we studied the popularity rankings of topics on 
2ch.net, one of the largest collected posting web pages in Japan,
and the book ranking of the amazon.co.jp, 
the Japanese counterpart of amazon.com, which is 
a large online bookstore 
quoted as one of the pioneering `long-tail' business in the era of
internet retails \cite{longtail}.
We performed a statistical fit of our model to the actual data, and
showed that we can apply to these social and economical activities 
the stochastic ranking process
with the (generalized) Pareto distribution as $\lambda$.
Statistical fits have shown \cite{HH072,HH073} that
these social and economical
activities are more `smash-hit' based rather than long-tail,
in contrast to the idea in \cite{longtail}.
The values of the Pareto parameter $0<b<1$ have also been found 
in a study of document access
in the MSNBC commercial news web sites \cite{Qiu} by directly counting 
the number of accesses.

Returning to the studies in MTF rules,
among the earliest works are \cite{Tsetlin1963,mv2frnt2,Letac1974},
where the formula for the stationary distributions of the MTF Markov chain 
is given.
Another earliest studies deals with the {\it search cost},
which is the 
position of the requested item before being moved.
(Figuratively, we can imagine a heap of reference papers.  
Every time we need a paper we start 
our search from the top of the heap and after use we return it on the top.) 
The formula of the average search cost for the stationary distribution
is first derived in \cite{mv2frnt1}.
Comparison of search cost probability with optimal ordering 
in the $N\to\infty$ limit is considered \cite{Jelenkovic99}.
The average search cost for stationary distribution has been studied 
in \cite{mv2frnt1,mv2frnt3} and
the comparison to that for the optimal ordering is found in 
\cite{mv2frnt3,Kingman75,Rivest1976,CHS88}.
A formula for generating function of the search cost 
is obtained in \cite{Fill96TCS}.
Search costs for non-stationary cases have also been studied 
\cite{Bitner79,Rodrigues95a,Fill96JTP,Fill96TCS}.
There are also studies of the 
conditional expectations of search costs \cite{mv2frnt4}, 
cache miss (fault) probability in the least-recently-used (LRU) caching
\cite{Fagin77,Jelenkovic99,Jelenkovic03,breslau,SM2006},
and the cases of generalized Zipf law or Pareto distribution 
as the jump-rate distribution 
\cite{Fill96TCS,Jelenkovic99,Jelenkovic03,SM2006,breslau}.
For summary of various studies of MTF models, see, for example, 
\cite{Fill96JTP,Jelenkovic99,SM2006}.

We will show in this paper that we can apply the mathematical results
in \cite{HH071,HH072} to derive formula for the asymptotic distribution 
of search cost $C_N$\,, for general jump-rate distribution $\lambda$.
A basic formula in the case of stationary distribution
is \eqnu{searchcostlimitstationarydistribution}:
\[
\limf{N} \prbseq{\infty}{\frac1NC_N> x} 
=\frac{\dsp \int_0^{\infty} e^{-wt_0(x)} w \lambda(dw)}{\dsp
 \int_0^{\infty} w\,\lambda(dw)}\,.
\]
Using the formula above, we can obtain the 
asymptotics of the search cost probabilities,
for general $\lambda$.
We have formula for non-stationary cases as well as the case of 
the stationary distribution (see
\eqnu{searchcostlimittransientdistribution}).

The plan of the present paper is as follows.
In \secu{2} we summarize the main mathematical results in \cite{HH071,HH072}.
In \secu{3} we use these results to derive the formula for 
the asymptotic distribution 
of search cost for general jump-rate distributions,
both for stationary and non-stationary cases.
In \secu{4} we reproduce and extend the formulas on asymptotics of the
search cost probabilities in the literature, using the results in \secu{3},
to show that our formula gives a unified way of deriving the results for 
the search costs in the MTF model.

\smallskip\par\textbf{Acknowledgements.}

The authors would like to thank Dr.~N.~Sugimine for bringing our
attention to the keyword, move-to-front rules.
The research of T.~Hattori 
is supported in part by KAKENHI 17340022 from 
the Ministry of Education, Culture, Sports, Science and Technology (MEXT).

\section{Stochastic ranking process.}
\seca{2}


Let $N$ be the total number of particles aligned in a queue
(records of information in a serial file, 
in terms of \cite{mv2frnt1},
or books on a single shelf, in terms of \cite{mv2frnt2,mv2frnt3}),
and for $i=1,2,\cdots,N$, and $t\ge0$,
let $X\ssN_i(t)$ be the position (ranking, in terms of 
\cite{HH071,HH072,HH073})
of particle $i$ in the queue at time $t$.

The particles jump at random jump times to the top position of the queue.
Denote by $\tau\ssN_{i,j}$, 
the time that particle $i$ jumps for the $j$-th time
to the top position.
Namely, for each $i$,
$X\ssN_i(\tau\ssN_{i,j})=1$, $j=1,2,\cdots$.
($\tau\ssN_{i,j}$ is the time of $j$-th request of record $i$, 
in terms of \cite{mv2frnt1}, or the time that a book is requested and 
returned at the left end of the shelf `nearest to the librarian's desk', 
in terms of \cite{mv2frnt2,mv2frnt3}.)
Besides the jump to the top position,
$X\ssN_i(t)$ changes its value when
some other particle nearer to the tail position jumps to the top
and the particle $i$ is pushed towards the tail 
to make room for the jumped particle.
Namely, for each $i'\ne i$ and $j'=1,2,\cdots$,
if $X\ssN_i(\tau\ssN_{i',j'}-0)<X\ssN_{i'}(\tau\ssN_{i',j'}-0)$ then
$\dsp X\ssN_i(\tau\ssN_{i',j'})=X\ssN_i(\tau\ssN_{i',j'}-0)+1$.
Otherwise, $X\ssN_i(t)$ is constant in $t$.

We assume that the jump times $\tau\ssN_{i,j}$ are independent in $i$,
and are independent of $X\ssN_i(t)$, $i=1,2, \cdots ,N$, $t\geq 0$.  
For simplicity 
of notation, we put
$\tau\ssN_{i,0}=0$, $i=1,2,\cdots,N$, and further assume that
for each $i=1,2,\cdots,N$,
$\{\tau\ssN_{i,j+1}-\tau\ssN_{i,j}\mid j=0,1,2,\cdots \}$
are independent whose distribution are identical for all $j$ and are
the exponential distribution
\eqnb
\eqna{stoppingtimeexponential}
\prb{\tau\ssN_i\le t} =1-e^{-w\ssN_it},\ t\ge 0\,,
\eqne
for a positive constant (the jump rate of the particle $i$) $w\ssN_i>0$.

Alternatively, we may define $X\ssN=(X\ssN_1,\cdots,X\ssN_N)$
as a Markov process on the state space $S_N$ of the set of $N!$
permutations of $\{1,2,\ldots ,N\}$,
with the Poisson jump times 
$\{\tau\ssN_{i,j}\mid i=1,2,\cdots,N,j=1,2,3,\cdots\}$
determined by \eqnu{stoppingtimeexponential}.

Note that with probability $1$,
$\tau\ssN_{i,j}$, $j=0,1,2,\cdots$, 
in \eqnu{stoppingtimeexponential}
is strictly increasing,
and that 
$\tau\ssN_{i,j} \ne \tau\ssN_{i',j'}$
for any different pair of suffices $(i,j)\ne (i',j')$, unless $j=j'=0$.
We may (and will) therefore work on the event that these properties
on $\tau\ssN_{i,j}$'s hold.
In particular,
if we align the distinct random times $\tau\ssN_{i,j}$
in an increasing order and denote the $k$-th number by $\sigma\ssN(k)$,
namely,
\eqnb
\eqna{jumptimes}
\arrb{l}\dsp
\{\sigma\ssN(k)\mid k=0,1,2,3,\cdots\}
=\{0\}\cup\{\tau\ssN_{i,j}\mid j=1,2,\cdots,\ i=1,2,\cdots,N\};
\\ \dsp
0=\sigma\ssN(0)<\sigma\ssN(1)<\sigma\ssN(2)<\cdots,
\arre
\eqne
then the stochastic chain
$Z\ssN(k)=(X\ssN_1(\sigma\ssN(k)),\cdots,X\ssN_N(\sigma\ssN(k))$, 
$k=0,1,2,\cdots$,
is a Markov chain on the state space of the permutations of $(1,2,\cdots,N)$,
satisfying the move-to-front rules of \cite{Tsetlin1963,mv2frnt1},
with the request probability $p\ssN_i$ of the record (or book) $i$ given by
\eqnb
\eqna{srp2mv2frntjumprate}
p\ssN_i=\prb{\sigma\ssN(1)=\tau\ssN_{i,1}}
=\int_0^{\infty} \prod_{j\ne i} e^{-w\ssN_jt}\,w\ssN_ie^{-w\ssN_it}dt
=\frac{w\ssN_i}{w\ssN_1+\cdots+w\ssN_N}\,.
\eqne
Note also that 
$\sigma\ssN(k+1)-\sigma\ssN(k)$, $k=1,2,\cdots$, are exponentially 
identically distributed
independent random variables, with a common distribution
\eqnb
\eqna{srptotaljumpdistr}
\prb{\sigma\ssN(1)\le t} =1-e^{-(w\ssN_1+\cdots+w\ssN_N)t},\ t\ge 0\,.
\eqne

Let, as in \cite{HH071},
$\dsp x\ssN_C(t)= \sharp\{ i\in \{1,2,\cdots,N\} \mid \tau\ssN_i \le t\}$
denote the boundary position in the queue
such that $\tau\ssN_i\le t$ if $X\ssN_i(t)\le x\ssN_C(t)$
and $\tau\ssN_i> t$ if $X\ssN_i(t)> x\ssN_C(t)$.
Namely, the particles towards the top side of $x\ssN_C(t)$ 
have experienced a jump by time $t$,
while none of the particles on the tail side of $x\ssN_C(t)$ has jumped
up to time $t$.

Denote the empirical distribution of jump rates by 
\eqnb
\eqna{lambdaN}
\lambda\ssN:=\frac1N \sum_{i=1}^N \delta_{w\ssN_i}\,,
\eqne
where, here and in the following, $\delta_c$ denotes a unit distribution
concentrated at $c$. Namely, for any set $A$,
\[
\int_A \delta_c (dw)  =\left\{ \arrb{ll}\dsp
1\,, & \mbox{if }\ \ c \in A, \\
\dsp
0\,, & \mbox{if }\ \ c \not\in A.
\arre \right. 
\]
\prpb[{\cite[Proposition 2]{HH071}}]
\prpa{yCt}
Assume
\eqnb
\eqna{HDLlimitPareto}
\lambda\ssN \to \lambda, \ N\to\infty,
\eqne
for a probability distribution $\lambda$ on $[0,\infty)$.
Then for $t\ge0$,
\eqnb
\eqna{yNCt}
y\ssN_C(t):=\frac1N x\ssN_C(t)
=\frac1N \sharp\{ i\in (1,2,\cdots,N) \mid \tau\ssN_i\le t\}
\eqne
converges in probability as $N\to\infty$ to
\eqnb
\eqna{yCt}
y_C(t)=1-\int_0^{\infty} e^{-wt} \lambda(dw).
\eqne
\DDD\prpe
This result says that the trajectory of a particle 
starting at the top position is
approximately given, for large $N$, by a deterministic trajectory
(adjusting the origin of the time parameter $t=0$ to be the time
that the particle is at the top position)
\eqnb
\eqna{NyCt}
Ny_C(t) =N(1-\int_0^{\infty} e^{-wt} \lambda(dw))
\sim N(1-\int_0^{\infty} e^{-wt} \lambda\ssN(dw))
= \sum_{i=1}^N (1-e^{-w\ssN_it}),
\eqne
as long as it remains in the queue (i.e., conditioned that it does not jump).
This is easy to recognize by noting that the motion of a particle in the
queue is caused by the random jumps of other particles, and that 
the law of large numbers replaces random jump times by their expectations.

We hereafter assume \eqnu{HDLlimitPareto}, together with
\eqnb
\eqna{lambda00}
\lambda(\{0\})=0,
\eqne
and
\eqnb
\eqna{lambdafiniteaverage}
\int_0^{\infty} w \lambda(dw)<\infty.
\eqne
As noted in \cite[Proposition 3]{HH071},
$y_C:\ [0,\infty)\to[0,1)$
then is continuous, strictly increasing, and bijective,
hence the inverse function  
$t_0:\ [0,1)\to[0,\infty)$ exists,
satisfying
\eqnb
\eqna{t0y}
y_C(t_0(y))=y, \ \ 0\le y<1\,,
\eqne
and
\eqnb
\eqna{yCasevaporatedparticlest0}
y=1-\int_0^{\infty} e^{-wt_0(y)} \lambda(dw).
\eqne
Differentiating \eqnu{yCt} and \eqnu{t0y}, we have
\eqnb
\eqna{dyCdt}
\diff{y_C}{t}(t)= \int_0^{\infty} we^{-wt} \lambda(dw)
=\frac{1}{\dsp \diff{t_0}{y}(y_C(t))}  \,.
\eqne

Now, consider an $N\to\infty$ scaling limit of the empirical distribution
on the product space of jump rate and position;
\eqnb
\eqna{HDLtN}
\eqna{empiricaldistribution}
\mu\ssN_{t}:=\frac1N \sum_i \delta_{(w\ssN_i,Y\ssN_i(t))}
\eqne
where,
\eqnb
\eqna{defprocj}
Y\ssN_i(t)=\frac1N\, (X\ssN_i(t)-1).
\eqne
We assume that the initial configuration of the queue
$(X\ssN_1(0),\cdots,X\ssN_N(0))=(x\ssN_{1,0},\cdots,x\ssN_{N,0})$
is such that the initial empirical distribution
$\mu\ssN_{0}$
converges weakly as $N\to\infty$ to a probability distribution $\mu_0$
whose second marginal is the Lebesgue measure on $[0,1)$;
for almost all $y\in[0,1)$, there exists a probability measure $\mu_{y,0}$
on the space of jump rates such that $\mu_0(dw,dy)=\mu_{y,0}(dw)\,dy$.

To state our main result in \cite{HH071},
We generalize \eqnu{yCt} and define
\eqnb
\eqna{ytildeCyt}
y_C(y,t)=1-\int_y^1 \int_0^{\infty} e^{-wt} \mu_{z,0}(dw)\,dz,
\ \ t\ge0,\ 0\le y< 1\,.
\eqne
In particular, $y_C(t)={y}_C(0,t)$.
For each $t\ge 0$,
${y}_C(\cdot,t):\ [0,1)\to[y_C(t),1)$ is a continuous, 
strictly increasing, bijective function of $y$,
hence the inverse function 
$\hat{y}(\cdot,t):\ [y_C(t),1)\to[0,1)$ exists:
\eqnb
\eqna{yhatyt}
1-y=\int_{\hat{y}(y,t)}^1 \int_0^{\infty} e^{-wt} \mu_{z,0}(dw)\,dz,
\ \ t\ge0,\ y_C(t)\le y< 1.
\eqne
In an analogy to \eqnu{NyCt},
the particle initially at the position $N y$,
will be approximately at $N y_C(y,t)$ at time $t$
for large $N$, provided the particle does not jump
to the top position by the time $t$.
It holds that
\eqnb
\eqna{hatyy}
\pderiv{\hat{y}}{y}(y,t)=\frac{1}{\dsp
 \int_0^{\infty} e^{-wt} \mu_{\hat{y}(y,t),0}(dw)}\,.
\eqne
\thmb[{\cite[Theorem 5]{HH071}}]
\thma{HDL}
Assume
\eqnu{HDLlimitPareto}, \eqnu{lambda00}, and \eqnu{lambdafiniteaverage}, and 
the convergence of the initial distribution $\mu\ssN_{0}$ as $N\to\infty$.
Then the joint empirical distribution $\mu\ssN_{t}(dw,dy)$
of jump rate and position at time $t$ 
converges as $N\to\infty$ to a distribution 
$\mu_{t}(dw,dy)=\mu_{y,t}(dw)\,dy$ on $\preals\times [0,1)$,
that is,
for any bounded continuous function 
$f:\ \preals\times[0,1) \to \reals$
\eqnb
\eqna{HDL}
\arrb{l}\dsp
\limf{N} \frac1N \sum_i f(w\ssN_i,Y\ssN_i(t) )
= \int_0^1 \left(\int_0^{\infty} f(w,y) \mu_{y,t}(dw)\right)\,dy,
\ \ \mbox{ in probability. }
\arre \eqne

The measure $\mu_{y,t}(dw)$ is given by
\eqnb
\eqna{Tets20070726}
\mu_{y,t}(dw)=\left\{ \arrb{ll}\dsp
\frac{\dsp
we^{-wt_0(y)} \lambda(dw)
}{\dsp
\int_0^{\infty} \tilde{w}e^{-\tilde{w}t_0(y)} \lambda(d\tilde{w})
}\,, & y<y_C(t), \\ \dsp
\frac{\dsp
e^{-wt} \mu_{\hat{y}(y,t),0}(dw)
}{\dsp
\int_0^{\infty} e^{-\tilde{w}t} \mu_{\hat{y}(y,t),0}(d\tilde{w})
}\,, & y>y_C(t).
\arre \right. 
\eqne
\DDD\thme
As noted in \cite[\S 2.1 Remark]{HH071},
the assumption \eqnu{lambdafiniteaverage}
assures that $\mu_{0,t}$ is well-defined.
The main results in \thmu{HDL} for $y>0$ hold
without \eqnu{lambdafiniteaverage}.

This completes a summary of main results in \cite{HH071}.

It is notationally simpler to write \eqnu{Tets20070726} in
a form integrated by $y$.
Recalling \eqnu{dyCdt} and \eqnu{hatyy}, we have
\eqnb
\eqna{Tets20070726integrated}
\mu_t(dw,[y,1))=\int_{z\in[y,1)}\mu_{z,t}(dw)\,dz
=\left\{ \arrb{ll} \dsp 
e^{-w t_0(y)} \lambda(dw), & y<y_C(t), \\ \dsp
\mu_0(dw,[\hat{y}(y,t),1)) e^{-wt}, & y>y_C(t).
\arre \right.
\eqne
Essential points about the formula are the importance of the curve $y=y_C(t)$,
and appearance of its inverse function $t_0$ as well as the inverse function 
$\hat{y}$ of $y_C(y,t)$.
An important observation in \cite{HH072} concerning these points is that
\eqnu{Tets20070726integrated} satisfies
a system of non-linear Burgers type partial differential equations
(see \eqnu{e1i} in \thmu{Burgers} below).
An initial value problem for \eqnu{e1i} is solved \cite{HH072} 
by a standard method
of characteristic curves, which precisely are the curves $y=y_C(t)$
and $y=y_C(y,t)$. The solution to the PDE is then written using
the inverse function of the characteristic curves.

To be explicit, consider, in particular, 
the case that the limit distribution of jump rates $\lambda$ is
a discrete distribution:
 $\dsp \lambda = \sum_{\alpha} \rho_{\alpha} \delta_{f_{\alpha}}$,
where the summation is taken over finite or countably infinite numbers,
or equivalently,
 $\dsp \lambda(\{f_{\alpha}\})= \rho_{\alpha}$, $\alpha=1,2,\cdots$,
where $\rho_{\alpha}$'s are positive numbers satisfying
 $\dsp \sum_{\alpha} \rho_{\alpha}=1$.
For $\alpha=1,2,\cdots$, put
\eqnb
\eqna{spacialdistributionfcn}
U_{\alpha}(y,t):=\mu_t(\{f_{\alpha}\},[y,1))
=\int_y^1 \mu_{z,t}(\{f_{\alpha}\})\,dz,
\eqne
and
$\dsp U_{\alpha}(y)=\int_y^1 \mu_{z,0}(\{f_{\alpha}\})\,dz$
for the initial data.
Then \eqnu{Tets20070726integrated} is written as
\eqnb
\eqna{Tets20070726discreteintegrated}
U_{\alpha}(y,t)=\left\{ \arrb{ll}\dsp
\rho_{\alpha}\, e^{-f_{\alpha} t_0(y)}
\,, & y<y_C(t), \\ \dsp
U_{\alpha}(\hat{y}(y,t))\, e^{-f_{\alpha}t} 
\,, & y>y_C(t).
\arre \right. 
\eqne
\thmb[{\cite[\S 2]{HH072}}]
\thma{Burgers}
Under the assumptions in \thmu{HDL},
\eqnu{Tets20070726discreteintegrated} is the unique (classical) solution 
to an initial value problem of a system of non-linear partial differential
equations defined by
\eqnb
\eqna{e1i}
\pderiv{U_{\alpha}}{t}(y,t) 
+\sum_{\beta} f_{\beta}\, U_{\beta}(y,t)\, \pderiv{U_{\alpha}}{y}(y,t)
=-f_{\alpha} U_{\alpha}(y,t),\ (y,t)\in[0,1)\times[0,\infty),
\ \alpha=1,2,\cdots,
\eqne
with the boundary condition
$\dsp U_{\alpha}(0,t)= \rho_{\alpha}\  \alpha=1,2,\cdots$, $t\ge 0$,
and the initial data 
$\dsp U_{\alpha}(y,0)=U_{\alpha}(y)$, $\alpha=1,2,\cdots$.
\DDD\thme
This completes a summary of the mathematical part of the main results 
in \cite{HH072}.

\section{Asymptotic distribution of search cost probabilities.}
\seca{3}

In this section,
we will relate our results summarized in \secu{2} to the previous studies
in move-to-front rules.

\subsection{Search cost.}

A typical quantity of interest in the studies of move-to-front rules
is the search cost $C_N$, which denotes the position of a particle
just before its jump to the top.

Let $Q\ssN_1$ be the random variable defined by
\eqnb
\eqna{Q1}
 \sigma\ssN(1)=\tau\ssN_{Q\ssN_1,1}
\eqne
where $\sigma\ssN$ is defined in \eqnu{jumptimes}.
Then $Q\ssN_1$ matches the definition of $Q_1$ in \cite{mv2frnt1},
and by definition,
\eqnb
\eqna{Q1distri}
\prb{Q\ssN_1=i}=p\ssN_i\,,\ i=1,\cdots,N,
\eqne
where $p\ssN_i$ is as in \eqnu{srp2mv2frntjumprate}.
$C_N$ (denoted by $X$ in \cite{mv2frnt1}) is then given by
$\dsp C_N=X\ssN_{Q\ssN_1}(\sigma\ssN(1)-0)$. Note that this is equal to
$\dsp X\ssN_{Q\ssN_1}(0)$,
because particles do not move before the first jump.
We see from \thmu{HDL} that, under the assumptions of \secu{2},
$C_N$ asymptotically scales as $N$ in the limit that $N\to\infty$,
and therefore the asymptotic properties of 
\eqnb
\eqna{searchcost}
 \frac1N C_N =Y\ssN_{Q\ssN_1}(0)
\eqne
where $Y\ssN_i$ is defined in \eqnu{defprocj},
is of interest.

\subsection{Distribution of search cost: Stationary case.}
\seca{3s}

As noted in \secu{2},
the stochastic ranking process can be viewed as a continuous-time Markov chain
on $S_N$.
Namely, $X\ssN(t)$ can be identified with an element 
$\pi =(\pi_1, \ldots ,\pi_N)$  of $S_N$ so that 
$\pi_i=X_i\ssN(t)$, $i=1, \ldots ,N$. 
The stochastic ranking process viewed as
a continuous-time Markov chain on $S_N$,
has the stationary distribution.
(The stationary distribution
is essentially the same as the stationary distribution 
of the move-to-front rules obtained by \cite{Tsetlin1963,mv2frnt2}
in a different way of correspondence,
$\pi _i$ being the label of the particle at the $i$-th position
in the references.)
Denote by $\EEseq{\infty}{}$ ($\prbseq{\infty}{}$, respectively) the
expectation (resp., probability) with respect to the stationary distribution
for the initial configurations.
If the distribution of the initial configuration
$(x\ssN_{1,0},\cdots,x\ssN_{N,0})=(X\ssN_1(0),\cdots,X\ssN_N(0))$
is the stationary distribution, then it is the distribution of 
$(X\ssN_1(t), \ldots ,X\ssN_N(t))$ for all $t\geq 0$.
In particular, 
for the $\mu\ssN_t$ in \eqnu{empiricaldistribution},
\eqnb
\eqna{stationaryempirical}
\mu\ssN_{\infty}:=
\EEseq{\infty}{\mu\ssN_0}=\EEseq{\infty}{\mu\ssN_t},\  \ t\ge0.
\eqne
Let $f(w,y)$ be a bounded continuous function with compact support.
Let $0<y_0<1$ be such that $f(w,y)=0$ for $y\ge y_0$\,,
and let $t>t_0(y_0)$, where $t_0$ is as in \eqnu{t0y}.
Note that $\mu_{y,t}$ in \eqnu{Tets20070726} for $t>t_0(y)$ is
constant in $t$ and independent of the initial distribution.
\thmu{HDL}, together with Fubini's Theorem and dominated convergence Theorem,
therefore implies
\eqnb
\eqna{stationarydistribution}
\arrb{l}\dsp
\limf{N} \int\!\!\!\int_{(w,y)\in[0,\infty)\times[0,1)}
 f(w,y)\, \mu\ssN_{\infty}(dw,dy)
=
\limf{N} \int\!\!\!\int_{(w,y)\in[0,\infty)\times[0,1)}
 f(w,y)\, \EEseq{\infty}{\mu\ssN_t(dw,dy)}
\\ \dsp {}
=
\EEseq{\infty}{  \int\!\!\!\int_{(w,y)\in[0,\infty)\times[0,1)}
 f(w,y)\, \mu_t(dw,dy) }
=
\frac{\dsp \int\!\!\!\int_{(w,y)\in[0,\infty)\times[0,1)}
 f(w,y) w e^{-wt_0(y)} dy \lambda(dw)}{\dsp
 \int_0^{\infty} w e^{-wt_0(y)}\,\lambda(dw)}
\arre \eqne
This implies that the joint empirical distribution $\mu\ssN_{\infty}$
of the jump rate and the position under the stationary distribution
in \eqnu{stationaryempirical} converges as $N\to\infty$ to
\eqnb
\eqna{Tets20070726stationary}
\limf{N} \mu\ssN_{\infty}(dw,dy)
=\mu_{\infty}(dw,dy):
= \frac{\dsp
we^{-wt_0(y)}dy \lambda(dw)
}{\dsp
\int_0^{\infty} \tilde{w}e^{-\tilde{w}t_0(y)} \lambda(d\tilde{w})}\,.
\eqne

The  distribution function of $\dsp \frac1NC_N$ in \eqnu{searchcost}
in the stationary state is then given by
\eqnb
\eqna{searchcostrandomdistribution}
\eqna{searchcostempiricaldistribution}
\arrb{l}\dsp
\prbseq{\infty}{\frac1NC_N> x} 
=\sum_{i=1}^N \prbseq{\infty}{Y_i\ssN(0)> x,\ Q\ssN_1=i}
\\ \dsp {}
=
\sum_{i=1}^N 
 \prbseq{\infty}{Q\ssN_1=i} \prbseq{\infty}{Y\ssN_i(0)> x}
=
\sum_{i=1}^N p\ssN_i \prbseq{\infty}{Y\ssN_i(0)> x}
\\ \dsp {}
= 
\sum_{i=1}^N
 w\ssN_i \frac{\dsp \prbseq{\infty}{Y\ssN_i(0)> x}}{\dsp \sum_{j=1}^N w\ssN_j}
=
\frac{\dsp \int\!\!\int_{(w,y)\in[0,\infty)\times(x,1)}
 w \,\mu\ssN_{\infty}(dw,dy)}{\dsp
 \int_0^{\infty} w\,\lambda\ssN(dw)},
\arre
\eqne
where, we first classified the total event by the first particle to jump, 
and then used the independence of $Q\ssN_1$ and $\{Y_i\ssN (0)\}$,
and finally, \eqnu{empiricaldistribution} and \eqnu{stationaryempirical}.
Combining \eqnu{Tets20070726stationary} with
\eqnu{searchcostrandomdistribution}, and changing the integration
variable $y$ to $t=t_0(y)$, using \eqnu{dyCdt}, we have
\eqnb
\eqna{searchcostlimitstationarydistribution1}
\eqna{searchcostlimitstationarydistribution}
\arrb{l} \dsp
\limf{N} \prbseq{\infty}{\frac1NC_N> x} 
=
\frac{\dsp \int\!\!\int_{(w,y)\in[0,\infty)\times(x,1)}
 w \,\mu_{\infty}(dw,dy)}{\dsp
 \int_0^{\infty} w\,\lambda(dw)}
\\ \dsp
=\frac{\dsp \int\!\!\int_{(w,t)\in[0,\infty)\times(t_0(x),\infty)}
 e^{-wt} dt\, w^2 \lambda(dw)}{\dsp
 \int_0^{\infty} w\,\lambda(dw)}
=\frac{\dsp \int_0^{\infty} e^{-wt_0(x)} w \lambda(dw)}{\dsp
 \int_0^{\infty} w\,\lambda(dw)}\,.
\arre \eqne
Similarly, we have, for a measurable function $f$,
\eqnb
\eqna{searchcostlimitstationarydistributionexpectation}
\arrb{l}\dsp
\limf{N} \EEseq{\infty}{f(\frac1N C_N)}
\\ \dsp {}
=\frac{\dsp \int\!\!\int_{(w,y)\in[0,\infty)\times[0,1)}
 wf(y) \,\mu_{\infty}(dw,dy)}{\dsp
 \int_0^{\infty} w\,\lambda(dw)}
=\frac{\dsp \int\!\!\int_{(w,t)\in[0,\infty)^2}
 f(y_C(t)) e^{-wt}\,dt\, w^2 \lambda(dw)}{\dsp
 \int_0^{\infty} w\,\lambda(dw)}\,.
\arre
\eqne
Note that if \eqnu{lambdafiniteaverage} fails,
then the denominator in the right hand side of
\eqnu{searchcostlimitstationarydistribution}
and \eqnu{searchcostlimitstationarydistributionexpectation}
diverges.

\subsection{Search cost: Comparison with optimally ordered case.}
\seca{RN}

Comparison between the search cost $C_N$ for the move-to-front rules and 
the search cost $R_N$ when the particles are in the optimal static ordering, 
i.e., when the particles are arranged in decreasing order of request 
probabilities $p_i$,
has been extensively studied \cite{mv2frnt3,CHS88,Jelenkovic99}.

For $0\le x\le 1$, define $w\ssN(x)$ by
\eqnb
\eqna{wxN}
\lambda\ssN([0,w\ssN(x)])=\frac{1}{N}[N(1-x)],
\eqne
where $[N(1-x)]$ denotes the largest integer not exceeding $N(1-x)$.
Noting \eqnu{Q1distri}, we have
\eqnb
\eqna{searchcostdistriorderedN}
\prb{\frac1N R_N>x}=\frac{\dsp \int_0^{w\ssN(x)} w\lambda\ssN(dw)}{\dsp 
\int_0^{\infty} w\lambda\ssN(dw)}\,.
\eqne
Taking ratio to \eqnu{searchcostempiricaldistribution},
and proceeding as in the derivation of 
\eqnu{searchcostlimitstationarydistribution},
we have
\eqnb
\eqna{CRtaildistribution}
\limf{N} \frac{\dsp\prbseq{\infty}{\frac1N C_N>x}}{\dsp\prb{\frac1N R_N>x}}
=\frac{\dsp \int_0^{\infty} e^{-wt_0(x)} w \lambda(dw)}{\dsp
 \int_0^{w(x)} w\lambda(dw)}\,,
\ 0< x< 1.
\eqne
where,
\eqnb
\eqna{wx}
\lambda([0,w(x)])=1-x.
\eqne

Note that all the $N\to\infty$ limit results so far, \textit{except for}
\eqnu{CRtaildistribution}, assume the condition
\eqnu{lambdafiniteaverage}, whereas 
\eqnu{CRtaildistribution} holds even if
\eqnu{lambdafiniteaverage} fails;
$\dsp \int_0^{\infty} w\, \lambda(dw)=\infty$.
(See the remark after \thmu{HDL}.)
Furthermore, if \eqnu{lambdafiniteaverage} holds,
then \eqnu{CRtaildistribution},
with \eqnu{yCasevaporatedparticlest0}, \eqnu{wx} and the dominated convergence 
theorem, implies 
\eqnb
\eqna{CRtaildistributionx0}
\lim_{x\to+0}
\limf{N} \frac{\dsp\prbseq{\infty}{\frac1N C_N>x}}{\dsp\prb{\frac1N R_N>x}}
=\frac{\dsp \int_0^{\infty} w \lambda(dw)}{\dsp \int_0^{\infty} w\lambda(dw)}
=1,
\eqne
which, considering a trivial equality
$\dsp
\prbseq{\infty}{\frac1N C_N\ge 0}=\prb{\frac1N R_N\ge 0}=1,
$
is a natural result.
In contrast, \eqnu{CRtaildistributionx0} may fail if 
$\dsp \int_0^{\infty} w\, \lambda(dw)=\infty$.
(See \secu{Pareto}.)

\subsection{Distribution of search cost: Non-stationary case.}
\seca{3ns}

We can generalize \eqnu{searchcostlimitstationarydistribution1}
in \secu{3s} to the non-stationary cases.
Let us return to the setting in \secu{2} and
assume that the initial value of the process is given:
$(X\ssN_1(0),\cdots,X\ssN_N(0))=(x\ssN_{1,0},\cdots,x\ssN_{N,0})$.
Let $\tau\ssN(t)=\inf \{\sigma\ssN(k)\mid \sigma\ssN(k) > t\}$
and define $I\ssN(t)$ by $\tau\ssN(t)=\tau\ssN _{I\ssN(t),j}$ for some $j$.
Define the search cost at time $t$ by $C_N(t)=X_{I\ssN(t)}(t)$.
We have, 
\eqnb
\eqna{searchcostgeneraldistribution}
\arrb{l}\dsp
\prbseq{t}{\frac1NC_N(t)> x} 
=\sum_{i=1}^N \prbseq{t}{Y_i(t)> x,\ 
I\ssN(t)=i}
=
\sum_{i=1}^N 
\prbseq{t}{Y\ssN_i(t)> x}\,  \prbseq{t}{I\ssN(t)=i}
\\ \dsp {}
= 
\sum_{i=1}^N
 w\ssN_i \frac{\dsp \prbseq{t}{Y\ssN_i(t)> x}}{\dsp \sum_{j=1}^N w\ssN_j}
=
\frac{\dsp \int\!\!\int_{(w,y)\in[0,\infty)\times(x,1)}
 w \,\mu\ssN_{t}(dw,dy)}{\dsp
 \int_0^{\infty} w\,\lambda\ssN(dw)}.
\arre
\eqne
Letting $N \to \infty $, we have
\eqnb
\eqna{searchcostlimittransientdistribution}
\limf{N} \prbseq{t}{\frac1NC_N(t)> x} 
=
\frac{\dsp \int\!\!\int_{(w,y)\in[0,\infty)\times(x,1)}
 w \,\mu_{y,t}(dw)\,dy}{\dsp
 \int_0^{\infty} w\,\lambda(dw)},
\eqne
where $\mu_{y,t}$ is given by \eqnu{Tets20070726}.

We also remark that since \eqnu{Tets20070726} coincides,
for $y<y_C(t)$ , with the stationary distribution
\eqnu{Tets20070726stationary},
the speed of approach to stationary state
is evaluated by \eqnu{yCt}:
\eqnb
\eqna{relaxation}
1-y_C(t)=\int_0^{\infty} e^{-wt} \lambda(dw).
\eqne

\section{Formulas related to search cost probabilities in the move-to-front
rules.}
\seca{4}

Some formulas related to the search cost for the move-to-front rules
have simple forms, and naturally was found in the early studies.
In this section we will derive formulas
corresponding to some of such nice formulas,
in the formulation of \secu{2}.

\subsection{Average search cost.}

\subsubsection{Asymptotic formula for the average search cost.}

In \cite{mv2frnt1}, 
the average search cost under the stationary distribution 
$\EEseq{\infty}{C_N}$ 
(denoted by $\mu$ in \cite{mv2frnt3})
is derived. 
Using the results and notations in \secu{3} and \secu{2},
we can calculate the asymptotics of this quantity.
With \eqnu{searchcostlimitstationarydistributionexpectation} we have
\[
\limf{N} \EEseq{\infty}{\frac1N C_N}
=
\frac1{\dsp \int_0^{\infty} w\,\lambda(dw)}
  \int\!\!\int_{(w,t)\in[0,\infty)^2}
  y_C(t)\, w^2 e^{-wt} \lambda(dw) \, dt.
\]
Using \eqnu{yCt} and performing the integration with respect to $t_0$,
we obtain
\eqnb
\eqna{searchcostexpectation}
\arrb{l}\dsp
\limf{N} \EEseq{\infty}{\frac1N C_N}
= \frac1{\dsp \int_0^{\infty} w\,\lambda(dw)}
\left( \int_0^{\infty}  w \lambda(dw)
-
  \int_0^{\infty}\int_0^{\infty}
 \frac{w^2}{w+\tilde{w}} \lambda(dw) \, \lambda(d\tilde{w})
\right)
\\ \dsp
\phantom{\limf{N} \EEseq{\infty}{\frac1N C_N}}
= 
\frac1{\dsp \int_0^{\infty} w\,\lambda(dw)}
  \int_0^{\infty}\int_0^{\infty}
 \frac{w\tilde{w}}{w+\tilde{w}} \lambda(dw) \, \lambda(d\tilde{w}).
 \arre
\eqne

Let us check that \eqnu{searchcostexpectation} is consistent with 
the corresponding result in \cite{mv2frnt1} (with notation changed to
those we adopt here):
\[
\EEseq{\infty}{C_N}=\frac12 + \sum_{i=1}^N \sum_{j=1}^N
 \frac{p\ssN_ip\ssN_j}{p\ssN_i+p\ssN_j}.
\]
With \eqnu{srp2mv2frntjumprate} and \eqnu{HDLlimitPareto} we have
\[ \arrb{l}\dsp
\EEseq{\infty}{\frac1N C_N}=\frac1{2N} + \frac1N \sum_{i=1}^N \sum_{j=1}^N
\frac{w\ssN_iw\ssN_j}{(w\ssN_i+w\ssN_j)(w\ssN_1+\cdots+w\ssN_N)}
\\ \dsp \phantom{\EE{\frac1N C_N}}
= \frac1{2N} + 
\frac1{\dsp \int_0^{\infty} w \lambda\ssN(dw)}
 \int_0^{\infty}\int_0^{\infty}
 \frac{w \tilde{w}}{w + \tilde{w}} \lambda\ssN(dw) \lambda\ssN(d\tilde{w})
\\ \dsp \phantom{\EE{\frac1N C_N}}
\to 
\frac1{\dsp \int_0^{\infty} w \lambda(dw)}
 \int_0^{\infty}\int_0^{\infty}
 \frac{w \tilde{w}}{w + \tilde{w}} \lambda(dw) \lambda(d\tilde{w}),
\ \ N\to \infty,
\arre \]
which coincides with \eqnu{searchcostexpectation}.

\subsubsection{Comparison with search cost for the optimal ordering.}

One of the first studies on comparison of the search cost $C_N$ with 
the search cost $R_N$ for the optimal ordering introduced in \secu{RN}
is found in \cite{mv2frnt3},
which gives a following universal bound for the expectations:
\[
\EEseq{\infty}{R_N}\le \EEseq{\infty}{C_N} \leq 2\EEseq{\infty}{R_N}-1.
\]
Corresponding relations for $N\to\infty$ then is
\eqnb
\eqna{searchcostexpectationcomparison}
\limf{N} \frac1N\EEseq{\infty}{R_N}\le 
\limf{N} \frac1N\EEseq{\infty}{C_N}\le
 2\limf{N} \frac1N\EEseq{\infty}{R_N}.
\eqne

To see that this relation follows from the results in \secu{3},
first note that
\[ \EEseq{\infty}{R_N}=\sum_{i=1}^N i p\ssN_i
= \sum_{(i,j);\; p\ssN_i \le p\ssN_j}  p\ssN_i
= \frac12 \sum_{i=1}^N \sum_{j=1}^N \min\{p\ssN_i,p\ssN_j \} +\frac{1}{2}\,.\]
With \eqnu{srp2mv2frntjumprate} and \eqnu{HDLlimitPareto} we then have
\[ \EEseq{\infty}{\frac1N R_N}
= \frac1{\dsp 2 \int_0^{\infty} w \lambda\ssN(dw)}
  \int_0^{\infty}\int_0^{\infty}
  \min\{w ,\tilde{w} \} \lambda\ssN(dw)\,\lambda\ssN(d\tilde{w}) +\frac{1}{2N},
\]
hence
\eqnb
\eqna{searchcostexpectationordered}
\limf{N} \EEseq{\infty}{\frac1N R_N}
= \frac1{\dsp 2 \int_0^{\infty} w\,\lambda(dw)}
  \int_0^{\infty}\int_0^{\infty}
 \min\{w ,\tilde{w} \} \lambda(dw) \, \lambda(d\tilde{w}).
\eqne
\eqnu{searchcostexpectationcomparison} is now a simple consequence of
\eqnu{searchcostexpectationordered} and \eqnu{searchcostexpectation},
if one notes a simple inequality
\[ \frac12 \min\{x,y\}\le \frac{xy}{x+y}\le \min\{x,y\},\ \ x\ge0,\ y\ge0. \]

We also note that there is a result \cite{CHS88} which proves that
a Hilbert's inequality implies a stronger universal upper bound,
which implies for the present case,
\eqnb
\eqna{searchcostexpectationcomparisonoptupperbd}
\limf{N} \frac1N\EEseq{\infty}{C_N}\le
 \frac{\pi}{2}\limf{N} \frac1N\EEseq{\infty}{R_N}.
\eqne
In fact, as derived in \cite{CHS88} we have,
\[ \arrb{l}\dsp
\frac12 \int_0^{\infty}\int_0^{\infty}
 \min\{w,\tilde{w}\}\lambda(dw)\lambda(d\tilde{w})
=\int_{w\le \tilde{w}} w \lambda(dw) \lambda(d\tilde{w})
=\int_0^{\infty} w \lambda([w,\infty)) \lambda(dw)
\\ \dsp {}
= -\left. \frac{w}{2} \lambda([w,\infty))^2\right|^{\infty}_0 + 
\frac12\int_0^{\infty}\lambda([w,\infty))^2dw
 = \frac12\int_0^{\infty} \lambda([w,\infty))^2dw, \arre \]
and
\[ \arrb{l}\dsp
\int_0^{\infty}\int_0^{\infty}
 \frac{w\tilde{w}}{w+\tilde{w}} \lambda(dw)\lambda(d\tilde{w})
\\ \dsp {}
=\int_0^{\infty}
\left[-\frac{w\tilde{w}}{w+\tilde{w}}\lambda([w,\infty))
\right]^{w=\infty}_{w=0} \lambda(d\tilde{w})
+\int_0^{\infty}\int_0^{\infty}
\left(\frac{\tilde{w}}{w+\tilde{w}}\right)^2
 \lambda(d\tilde{w})\lambda([w,\infty))dw
\\ \dsp {}
= - \int_0^{\infty}
\biggl[\biggl(\frac{\tilde{w}}{w+\tilde{w}}\biggr)^2\lambda([\tilde{w},\infty))
\biggr]^{\tilde{w}=\infty}_{\tilde{w}=0} \lambda([w,\infty))dw
+\int_0^{\infty}\int_0^{\infty}
 \frac{2w\tilde{w}}{(w+\tilde{w})^3}
\lambda([w,\infty))\lambda([\tilde{w},\infty))dwd\tilde{w}
\\ \dsp {}
=\int_0^{\infty}\int_0^{\infty}
 \frac{2w\tilde{w}}{(w+\tilde{w})^3}
\lambda([w,\infty))\lambda([\tilde{w},\infty))dwd\tilde{w},
\arre \]
which, with the Hilbert's inequality in the form 
\cite[\S 9.3]{Inequalities}
for $\dsp K(x,y)=\frac{4xy}{(x+y)^3}$, $p=q=2$, and $g=f\ge 0$;
\[ \arrb{l}\dsp
\int_0^{\infty}\int_0^{\infty} \frac{4xy}{(x+y)^3} f(x)f(y) dxdy
 \le k \int_0^{\infty} f(x)^2dx;
\ \ k=\int_0^{\infty}K(x,1)\frac{dx}{\sqrt{x}}
=\frac{4\Gamma(\frac32)^2}{\Gamma(3)}=\frac{\pi}2,
\arre \]
imply
\eqnu{searchcostexpectationcomparisonoptupperbd}.

\subsubsection{Conditional expectations of search costs.}

In \cite{mv2frnt3},
the average search cost conditioned on specific particle $i$
(denoted by $\mu_i$ in the reference), has been obtained.  It is related to 
$\EEseq{\infty}{C_N}$ by 

\eqnb
\eqna{conditionalaveragesearchcost}
\EEseq{\infty}{C_N}=\sum_{i=1}^N p\ssN_i \mu_i\,.
\eqne
In terms of the conditional expectation
$\dsp\EEseq{\infty}{C_N\mid Q\ssN_1}$,
conditioned on the sigma algebra
\[
 \sigma[Q\ssN_1]=\sigma[\{\sigma\ssN(1)=\tau\ssN_{i,1}\},\ i=1,2,\cdots,N]
\]
(recall \eqnu{Q1}),
we have
\eqnb
\eqna{particlesearchcost}
\EEseq{\infty}{C_N\mid Q\ssN_1}(\omega)=\mu_i\,,\ \mbox{ if }
\ Q\ssN_1(\omega)=i.
\eqne
With \eqnu{srp2mv2frntjumprate} 
we reproduce \eqnu{conditionalaveragesearchcost}.

In considering such quantities,
we naturally come across the distribution of `jumped particles', that is, 
the distribution of $Q\ssN_1$\,.
Note that the time evolution of the system is 
dependent only on the jump rates.
Therefore the search cost of particle $i$ in the stationary state is 
dependent on $i$ only through its jump rate $w\ssN_i$;
if $w\ssN_i=w\ssN_j$ then the search cost for $i$ and $j$ has the same
distribution.
In particular,
\eqnb
\eqna{conditionalsearchcosts}
\EEseq{\infty}{C_N\mid Q\ssN_1}=\EEseq{\infty}{C_N\mid W_N},
\eqne
where $W_N=w\ssN_{Q\ssN_1}$. 

Proceeding as in the argument for 
\eqnu{searchcostrandomdistribution}, we have, for a bounded measurable 
function $f$,
\[ \arrb{l}\dsp
\EEseq{\infty}{f(W_N)}
=\sum_{i=1}^N  f(w\ssN_i) \prbseq{\infty}{Q\ssN_1=i}
\\ \dsp
=\sum_{i=1}^N p\ssN_i f(w\ssN_i)
=\frac{\dsp  \int\!\!\int_{(w,y)\in[0,\infty)\times[0,1)}
 f(w) w \,\mu\ssN_{\infty}(dw,dy)}{\dsp
 \int_0^{\infty} w\,\lambda\ssN(dw)}.
\arre \]
As in \eqnu{searchcostlimitstationarydistribution},
\thmu{HDL} therefore implies, for a bounded continuous function $f$
\eqnb
\eqna{jumpedparticle}
\limf{N} \EEseq{\infty}{f(W_N)}
=
\frac{\dsp  \int\!\!\int_{(w,y)\in[0,\infty)\times[0,1)}
 f(w) w \,\mu_{\infty}(dw,dy)}{\dsp
 \int_0^{\infty} w\,\lambda(dw)}
=
\frac{\dsp \int_0^{\infty} f(w) w \,\lambda(dw)}{\dsp
 \int_0^{\infty} w\,\lambda(dw)}\,.
\eqne
In other words, the distribution of the jumped particle jump rates
in the stationary state converges weakly to a probability measure
$\dsp \frac{\dsp w \,\lambda(dw)}{\dsp
 \int_0^{\infty} \tilde{w}\,\lambda(d\tilde{w})}$,
as $N\to\infty$.

Since $t_0(0)=0$, this distribution is equal to
$\mu_{0,\infty}$ in \eqnu{Tets20070726stationary},
which is the distribution at the top end of the queue.
An intuitive meaning of this equality is that
the jumped particles jump to the top position
(the requested records are placed at the top position)
so the distribution at $y=0$ is the distribution of the jumped particles.


As noted in \eqnu{conditionalsearchcosts},
to obtain the
average search cost of a specific particle $i$ 
(denoted by $\mu_i$ in \cite{mv2frnt3}),
it suffices to calculate the average search cost
conditioned on the jump rate of the jumped particle
$\dsp f(W_N)=\EEseq{\infty}{C_N\mid W_N}$.
A basic property of conditional expectation, with
\eqnu{searchcostexpectation} and \eqnu{jumpedparticle}, implies
\[\arrb{l}\dsp
 \frac1{\dsp \int_0^{\infty} w\,\lambda(dw)}
  \int_0^{\infty}\int_0^{\infty}
 \frac{w\tilde{w}}{w+\tilde{w}} \lambda(dw) \, \lambda(d\tilde{w})
=\limf{N} \EEseq{\infty}{\frac1N C_N}
=
\limf{N} \EEseq{\infty}{\EEseq{\infty}{\frac1N C_N\mid W_N}}
\\ \dsp {}
=
\frac1{\dsp \int_0^{\infty} w\,\lambda(dw)}\,
 \int_0^{\infty} \limf{N} \EEseq{\infty}{\dsp \frac1N C_N\mid W_N}(w)\, w
 \,\lambda(dw).
\arre \]
Thus we find
\eqnb
\eqna{searchcostconditionalaverage}
\limf{N} \frac1N \EEseq{\infty}{C_N\mid W_N}(w)\,
= \int_0^{\infty} \frac{\tilde{w}}{w+\tilde{w}}\, \lambda(d\tilde{w}).
\eqne

This result is to be compared with $\mu_i$ in \cite[Eq.~(10)]{mv2frnt3},
which reads in our notation,
\[
\frac1N \EEseq{\infty}{C_N\mid W_N}(w_i)
=\frac1N \mu_i
= \frac1{2N} +\frac1N\sum_{j=1}^N \frac{w\ssN_j}{w\ssN_i+w\ssN_j}\,.
\]
For large $N$, \eqnu{HDLlimitPareto} then implies
\[
\frac1N \EEseq{\infty}{C_N\mid W_N}(w\ssN_i)
= \frac1{2N}
 + \int_0^{\infty}
 \frac{\tilde{w}}{w\ssN_i+\tilde{w}} \, \lambda\ssN(d\tilde{w})
\sim
 \int_0^{\infty}
 \frac{\tilde{w}}{w_i+\tilde{w}} \, \lambda(d\tilde{w}),
\]
which is consistent with \eqnu{searchcostconditionalaverage}.

\subsection{Cache miss probability.}

If $x\le y_C(t)$ we can reduce
\eqnu{searchcostlimittransientdistribution}
further and have
\eqnb
\eqna{searchcostlimittransientdistributionxleyC}
\limf{N} \prbseq{t}{\frac1NC_N(t)> x}
=\frac{\dsp \int_0^{\infty} e^{-wt_0(x)} w \lambda(dw)}{\dsp
 \int_0^{\infty} w\,\lambda(dw)}\,.
\eqne
This is because the limiting distribution $\mu_{y,t}$ 
for $y<y_C(t)$ is equal to that for stationary case $\mu_{y,\infty}$.
(See \eqnu{Tets20070726} and \eqnu{Tets20070726stationary}.)
Hence we have, for $x\le y_C(t)$,
\[ \arrb{l} \dsp
\limf{N} \prbseq{t}{\frac1N C_N>x}=
1-\limf{N} \prbseq{t}{\frac1N C_N\le x}=
1-\limf{N} \prbseq{\infty}{\frac1N C_N\le x}
\\ \dsp {}
=\limf{N} \prbseq{\infty}{\frac1N C_N>x},
\arre \]
so that
\eqnu{searchcostlimitstationarydistribution} implies
\eqnu{searchcostlimittransientdistributionxleyC}.

The cache miss (fault) probability in the least-recently-used (LRU) caching
has been one of the modern area of extensive study in the application of 
the move-to-front rules 
\cite{Fagin77,Jelenkovic99,Jelenkovic03,breslau,SM2006}.
If there is $N$ records of information in a computer memory, or $N$ web pages
on the internet, out of which $Nx$ records or pages, respectively,
can be cached for a further quick access,
the event $C_N>Nx$ represents cache miss or cache fault,
by regarding particles as records of information or web pages to be accessed.
The probability
\eqnu{searchcostlimittransientdistributionxleyC}
is therefore of interest.

In particular, \cite{breslau} considers a quantity, defined, in our notation,
by
\eqnb
\eqna{missprobabN}
M\ssN(t)= \prbseq{t}{\frac1NC_N> y\ssN_C(t)}.
\eqne
Recalling the definition \eqnu{yNCt} of $y\ssN_C$,
we see that $M\ssN(t)$ is the probability that the jump at time $t$ 
is the jumped particle's first jump since $t=0$.
$M\ssN(t)$ therefore corresponds to the cache miss (fault) probability
in an ideal case that all the once requested records are stored in 
a cache memory of ideally large size.

Since the limiting distribution \eqnu{searchcostlimittransientdistribution}
of $\dsp\frac1N C_N$ is continuous and
and $y_C^{(N)}(t)$ converges in probability to $y_C(t)$, we have 
\eqnb
\eqna{searchcostlimittransientdistributionxeqCN}
\limf{N} M\ssN(t)
= \limf{N} \prbseq{t}{\frac1NC_N(t)> y_C(t)}.
\eqne
Substituting $x=y_C(t)$ in
\eqnu{searchcostlimittransientdistributionxleyC},
we have
\eqnb
\eqna{missprobability}
M(t) :=\limf{N} M\ssN(t)
=\limf{N} \prbseq{t}{\frac1NC_N(t)> y_C(t)}
=\frac{\dsp \int_0^{\infty} e^{-w t} w \lambda(dw)}{\dsp
 \int_0^{\infty} w\,\lambda(dw)}\,.
\eqne
Note that $M(t)$ is independent of the initial configuration $\mu\ssN_0$.

\subsection{Case of generalized Zipf law or Pareto distribution.}
\seca{Pareto}


In the preceding subsections, we dealt with formula for an arbitrary
distribution of the jump rates $\lambda$.
In the literature, there are formula for specific request probabilities,
among which the generalized Zipf law (also known as
power-law) is of importance in practical
applications.
Let $a$ and $b$ be positive constants and consider the jump rates
\eqnb
\eqna{Paretodiscrete}
w_i=a\left(\frac{N}{i}\right)^{1/b}, \ i=1,2,3,\cdots,N.
\eqne
In applying to move-to-front rules, $a=w_N$ is the smallest jump rate
and $\dsp b=\frac{\log N}{\dsp \log\frac{w_1}{w_N}}$ is an exponent
representing the equality of jump rates among the particles.

In \cite{HH072,HH073} we studied the rankings 
of 2ch.net and amazon.co.jp. 
2ch.net is one of the largest collected posting web pages in Japan.
Posting web pages are classified by categories (`boards'), 
and each category has a list of topics of posting web pages (`threads'). 
These lists are updated by the `last-written-thread-at-the-top" rule. 
Amazon.co.jp is the Japanese counterpart of amazon.com, which is
a large online bookstore 
quoted as one of the pioneering `long-tail' business in the era of
internet retails \cite{longtail}.  They show sales ranks of all the books on 
their catalogs.
We have shown that we can apply the stochastic ranking process
with the (generalized) Pareto distribution 
for the distribution of jump rates
in these social and economical activities, and
by performing statistical fits of the data from these web results,
we extracted the index $b$ in \eqnu{Paretodiscrete}.
We obtained $b=0.61$ for 2ch.net and $b=0.81$ for amazon.co.jp,
both indicating $0<b<1$, which implies that these social and economical
activities are more `smash-hit' based rather than long-tail,
in contrast to the idea in \cite{longtail}.
The values in $0<b<1$ has also been found in a study of document access
in the MSNBC commercial news web sites \cite{Qiu} by direct measurements
(that is, the distribution $\lambda$ is directly measurable in the study
of \cite{Qiu} and a theory of move-to-front rules is unnecessary).

Let us turn to the search cost probabilities.
$\lambda$ of \eqnu{HDLlimitPareto} is readily calculated:
\eqnb
\eqna{lambdaPareto}
\lambda([0,w]) = \left\{ \arrb{ll} 0,& 0\le w<a,\\
 \dsp 1-\left(\frac{a}{w}\right)^b,& w\ge a. \arre \right.
\eqne
The continuous distribution $\lambda$ determined by \eqnu{lambdaPareto}
is called the (generalized) Pareto distribution \cite{Pareto1896} 
(or log-linear distribution), especially in social studies, 
and is used to explain various social distributions, 
typically that of incomes.

With \eqnu{wx}
we have $\dsp w(x)= ax^{-1/b}$, and
the denominator in the right hand side of \eqnu{CRtaildistribution} is
\eqnb
\eqna{CRtaildistributionPareto1}
 \int_0^{w(x)} w\lambda(dw) = \frac{ab}{1-b} (x^{1-1/b}-1).
\eqne
For the numerator of \eqnu{CRtaildistribution} we have
\eqnb
\eqna{CRtaildistributionPareto2}
\int_0^{\infty} e^{-wt_0(x)} w \lambda(dw)
=\int_a^{\infty} e^{-wt_0(x)}\, b \left(\frac{a}{w}\right)^b\,dw
=\frac{b}{t_0(x)}\, (at_0(x))^b \, \Gamma(1-b,at_0(x)),
\eqne
where $\dsp \Gamma(z,p)=\int_{p}^{\infty} e^{-w} w^{z-1} \, dw$
is the incomplete Gamma function.
To evaluate this, we recall \eqnu{yCasevaporatedparticlest0}
and perform integration by parts, to find
\eqnb
\eqna{CRtaildistributionPareto3}
1-x=\int_0^{\infty} e^{-wt_0(x)}\frac{ba^b}{w^{b+1}}\,dw
=e^{-at_0(x)} - (at_0(x))^b \Gamma(1-b,at_0(x)).
\eqne
Substituting \eqnu{CRtaildistributionPareto1},
 \eqnu{CRtaildistributionPareto2},
and \eqnu{CRtaildistributionPareto3} in \eqnu{CRtaildistribution}, we have
\eqnb
\eqna{CRtaildistributionPareto}
\limf{N} \frac{\dsp\prbseq{\infty}{\frac1N C_N>x}}{\dsp\prb{\frac1N R_N>x}}
=\frac{1-b}{at_0(x)}\frac{\dsp e^{-at_0(x)}-1+x}{\dsp x^{1-1/b}-1}\,.
\eqne
This formula is valid for all $b>0$ and $0<x<1$.

Concerning the condition \eqnu{lambdafiniteaverage},
we see from \eqnu{lambdaPareto},
\eqnb
\eqna{lambdafiniteaveragePareto}
\int_0^{\infty} w\lambda(dw) = \frac{ab}{b-1} \,,
\eqne
so that \eqnu{lambdafiniteaverage} is equivalent to $b>1$ for the 
Pareto distribution.
Hence, as discussed in \secu{RN}, \eqnu{CRtaildistributionx0} holds if $b>1$.
In contrast, if $0<b<1$, then noting $\dsp \lim_{x\to+0} t_0(x)=0$
(which is seen from the definition \eqnu{yCasevaporatedparticlest0}),
we have
$\dsp \lim_{x\to+0} \Gamma(1-b,at_0(x))=\Gamma(1-b)$, and
\eqnu{CRtaildistributionPareto3} implies
\[  at_0(x) \sim \left(\frac{x}{\Gamma(1-b)}\right)^{1/b},\ \ x\to 0
\ \ \mbox{ if $0<b<1$}, \]
and \eqnu{CRtaildistributionPareto} then implies
\eqnb
\eqna{CRtaildistributionParetox0}
\lim_{x\to+0}
\limf{N} \frac{\dsp\prbseq{\infty}{\frac1N C_N>x}}{\dsp\prb{\frac1N R_N>x}}
=(1-b)\, \Gamma(1-b)^{1/b}.
\eqne
The quantity in the right hand side of this result is obtained 
in \cite[Theorem 3]{Jelenkovic99}.
Note that the reference formulates $N=\infty$ case from the beginning
(in our notation, this is attained by letting $a$ to be proportional
to $N^{-1/b}$ in \eqnu{Paretodiscrete}), and 
a limit $n\to \infty$ is taken in Theorem 3 of \cite{Jelenkovic99}.
We begin with $N\to\infty$, fixing $x$, and then take $x\to 0$ limit
in \eqnu{CRtaildistributionParetox0}.
Rigorously speaking, these are different limits and
\eqnu{CRtaildistributionParetox0} is a new result.
However, since our $x$ and $n$ in \cite{Jelenkovic99} are related by
$n=Nx$ when $N<\infty$, both results are consistently talking about
`large $N$, large $n$, and small $x$' for $0<b<1$.

Concerning \eqnu{searchcostlimitstationarydistributionexpectation},
a general formula for the expectation of search cost,
we have
\eqnb
\eqna{searchcostlimitstationarydistributionexpectationPareto}
\limf{N} \EEseq{\infty}{f(\frac1N C_N)}
= (b-1) (at)^{b-2} a \int_0^{\infty} f(y_C(t)) \Gamma(2-b,at) dt,
\eqne
where, \eqnu{yCt} implies
\eqnb
\eqna{yCtPareto}
y_C(t)=1-b(at)^b\Gamma(-b,at).
\eqne
Noting that
\[
\diff{y_C}{t}(t)=ab (at)^{b-1} \Gamma(1-b,at)
\]
and an integration by parts formula 
\[
\Gamma(z+1,p)=e^{-p}p^z + z\Gamma(z,p)
\]
for the incomplete gamma function,
we have another expression
\eqnb
\eqna{searchcostlimitstationarydistributionexpectationPareto2}
\limf{N} \EEseq{\infty}{f(\frac1N C_N)}
= 
(b-1)\int_0^{\infty} f(y_C(t)) e^{-at} \frac{dt}{t}
-\frac{(b-1)^2}{ab} \int_0^1 \frac{f(y)}{t_0(y)} dt.
\eqne
It seems, however, difficult to simplify the formula for general $f$.

Concerning the miss probability $M(t)$ of \eqnu{missprobability},
the denominator is finite if $b>1$,
and we have, after an integration by parts,
\[
M(t)= (b-1) (at)^{b-1} \int_{at}^{\infty} e^{-x} x^{-b} dx
= e^{-at} - (at)^{b-1} \Gamma(2-b,at).
\]
For $1<b<2$ this implies
\eqnb
\eqna{cachemissParetobgt1}
M(t) = 1- \Gamma(2-b) (at)^{b-1} +O(at),\  \ t\to 0.
\eqne

In \cite{breslau} the web caching is studied, in which
the hit-ratio for the $R$-th request is defined, in our notation, by
\[
H\ssN(R)=1-M\ssN(\sigma\ssN(R)),
\]
where $M\ssN(t)$ is defined in \eqnu{missprobabN} and
$\sigma\ssN(k)$ in \eqnu{jumptimes}.
With \eqnu{cachemissParetobgt1} and properties of $\sigma\ssN$
(see \eqnu{srptotaljumpdistr}), together with law of large numbers,
we see that
$\dsp H(R)=\limf{N} H\ssN(NR)$ scales as $R^{b-1}$.
This is consistent with the argument in \cite{breslau}
which claims $H(R) \propto R^{b-1}$ for $1\ll R\ll N$.
The reference further obtains $1/b=0.83-0.90$ ($b=1.11-1.20$)
using actual web data.

\end{document}